\documentclass[12pt]{article}
\usepackage{amsmath,amsthm,amsfonts,latexsym,amssymb,fancyhea,supertab}
\usepackage{a4}

  \def\Q{\mathbb{Q}}
\def\Z{\mathbb{Z}} \def\R{\mathbb{R}}

\def\cE{{\cal E}}

\def\al{\alpha} \def\be{\beta} \def\ep{\epsilon} 
\def\th{\theta}  \def\ga{\gamma}
\def\la{\lambda}

\begin{document}

\newtheorem*{u8_result}{Theorem}%[section]
\newtheorem{problem}{Problem}%[section]
\newtheorem{maintheorem}{Theorem}
\title{On squares in Lucas sequences}
\author{A.~Bremner\thanks{Department of Mathematics, Arizona State University, Tempe AZ, USA, e-mail: bremner@asu.edu, http://\~{}andrew/bremner.html} \and 
N.~Tzanakis\thanks{Department of Mathematics, University of Crete,
Iraklion, Greece, e-mail: tzanakis@math.uoc.gr , http://www.math.uoc.gr/\~{}tzanakis}
}
\date{July 5, 2006}
\maketitle

\begin{abstract}
\noindent
Let $P$ and $Q$ be non-zero integers. The Lucas sequence $\{U_n(P,Q)\}$
is defined by
$U_0=0, \quad U_1=1, \quad U_n= P U_{n-1}-Q U_{n-2} \quad (n \geq 2).$
The question of when $U_n(P,Q)$ can be a perfect square has generated interest in the literature.
We show that for $n=2,...,7$, $U_n$ is a square for infinitely many pairs $(P,Q)$ with
$\gcd(P,Q)=1$; further, for $n=8,...,12$, the only non-degenerate sequences where $
\gcd(P,Q)=1$
and $U_n(P,Q)=\Box$, are given by $U_8(1,-4)=21^2$, $U_8(4,-17)=620^2$, and $U_{12}(1,-1)=12^2$.
\end{abstract}
Keywords: Lucas sequence, squares, genus two curves

\section{Introduction} \label{introduction}
Let $P$ and $Q$ be non-zero integers. The Lucas sequence $\{U_n(P,Q)\}$ 
is defined by
\begin{equation}
\label{Lucas}
U_0=0, \quad U_1=1, \quad U_n= P U_{n-1}-Q U_{n-2} \quad (n \geq 2).
\end{equation}
The sequence $\{U_n(1,-1)\}$ is the familiar Fibonacci sequence, and it was proved by 
Cohn~\cite{Co1} in 1964 that the only perfect square
greater than $1$ in this sequence is $U_{12}=144$. The question arises, 
for which parameters $P$, $Q$, can $U_n(P,Q)$ be a perfect square? 
In what follows, we shall assume that we are not dealing with the degenerate 
sequences corresponding to $(P,Q)=(\pm 1,1)$, where $U_n$ is periodic
with period $3$, and we also assume $(P,Q) \neq (-2,1)$ (in which case 
$U_n=\Box$ precisely when $n$ is an odd square) and  $(P,Q) \neq (2,1)$ 
(when $U_n=\Box$ precisely when $n$ is square).
Ribenboim and McDaniel~\cite{RM1} with only elementary methods show that when $P$ and $Q$ are 
{\it odd}, and $P^2-4Q>0$,
then $U_n$ can be square only for $n=0,1,2,3,6$ or $12$; and that there are at most two indices 
greater than 1 for which $U_n$ can be square.
They characterize fully the instances when $U_n=\Box$, for $n=2,3,6$.
Bremner \& Tzanakis~\cite{BT} extend these results by determining all Lucas 
sequences $\{U_n(P,Q)\}$ with $U_{12}=\Box$, subject only to the restriction 
that $\gcd(P,Q)=1$ (it turns
out that the Fibonacci sequence provides the only example). Under the same hypothesis, all 
Lucas sequences with $\{U_n(P,Q)\}$ with $U_9=\Box$ are determined.
There seems little mention in the literature of when
under general hypotheses $U_n(P,Q)$ can be a perfect square. 
It is straightforward to see from Theorem 1 of Darmon and Granville~\cite{DG}
that for $n$ sufficiently large ($n \geq 11$ certainly suffices), 
the equation $U_n(P,Q)=\Box$ can have only finitely many
solutions for coprime $P,Q$.
Note that 
for $n \geq 1$, $U_n(k P, k^2 Q)=k^{n-1} U_n(P,Q)$, and so for fixed $P$, $Q$, 
and {\it even} $n$, appropriate choice of $k$ gives a sequence
with $U_n(k P, k^2 Q)$ a perfect square. The restriction to $\gcd(P,Q)=1$ is therefore
a sensible one, and we shall assume this from now on.  Rather curiously, 
a small computer search
reveals sequences with $U_n(P,Q)$ a perfect square only for 
$n=0,\ldots,8$, and $n=12$. In this paper, we shall dispose of this 
range of $n$.  Bremner \& Tzanakis~\cite{BT} have addressed the cases 
$n=9, 12$. Section \ref{up to U7}
of this paper addresses the case $U_n(P,Q)=\Box, \; n\leq 7$,
which can be treated entirely elementarily. Section \ref{U8 to U11}
addresses the cases $U_n(P,Q)=\Box$, $8 \leq n \leq 11$.
In these instances, we deduce a finite collection of curves, whose 
rational points cover all required solutions. In turn, this
reduces to a number of problems of similar type, namely, finding all 
points on an elliptic curve defined over a number field $K$ subject 
to a ``$\Q$-rationality" condition on the $X$-coordinate. 
Nils Bruin has powerful
techniques for addressing this type of problem, 
and~\cite{Br1},~\cite{Br2},~\cite{Br3},~\cite{Br4} provide details and examples. See 
in particular \S4 of~\cite{Br4} for development of the underlying mathematics.
The latest release of Magma~\cite{Mag} now contains Bruin's routines 
and so we only set up the appropriate computation here, with details of
the Magma programs available on request. \\ \\
The results of Sections \ref{up to U7}, \ref{U8 to U11}, when
combined with the results of Bremner \& Tzanakis~\cite{BT} give the following theorem:
\begin{maintheorem} 
Let $P,Q$ be non-zero coprime integers such that if $Q=1$ then 
$P \neq \pm 1,\pm2$ (that is, $P,Q$ determine a non-degenerate Lucas sequence).
Then (i) for $n=2,...,7$, $U_n(P,Q)$ is a square for infinitely many such 
pairs $(P,Q)$, and (ii) for $n=8,...,12$, the only solutions of 
$U_n(P,Q)=\Box$ are given by $U_8(1,-4)=21^2$, $U_8(4,-17)=620^2$, 
and $U_{12}(1,-1)=12^2$.
\end{maintheorem}

\section{Solution of $U_n(P,Q)=\Box, \; n\leq 7$} \label{up to U7}
Certainly $U_2(P,Q)=\Box$ if and only if $P=a^2$, and $U_3(P,Q)=\Box$ if and
only if $P^2-Q=a^2$.\\
Now $U_4(P,Q)=\Box$ if and only if $P(P^2-2Q)=\Box$, so if and only if either 
$P= \delta a^2, Q=\frac{1}{2}(a^4-\delta b^2)$, or $P=2\delta a^2, Q=2a^4-\delta b^2$,
with $\delta=\pm 1$ (where, in the first instance, $a b$ is odd and in the
second instance $b$ is odd). \\
The demand that $U_5(P,Q)$ be square is that $P^4-3 P^2 Q+Q^2 = \Box$, equivalently,
that $1-3 x+x^2=\Box$, where $x=Q/P^2$. Parametrizing the quadric,
$Q/P^2 = (5 \lambda^2+6 \lambda \mu +\mu^2)/(4 \lambda \mu)$, where, without loss
of generality, $(\lambda, \mu)=1$, $\lambda>0$, and $\mu \not \equiv 0 \pmod 5$. Necessarily
$(\lambda, \mu)=(a^2, \pm b^2)$, giving $(P,Q)=(2 a b, 5 a^4 + 6 a^2 b^2 +b^4)$ or
$(2 a b, -5 a^4+6 a^2 b^2 - b^4)$  if $a$ and $b$ are of opposite parity, and 
$(P,Q)=(a b, \frac{1}{4}(5 a^4 + 6 a^2 b^2 +b^4))$ or 
$(a b, \frac{1}{4}(-5 a^4 + 6 a^2 b^2 - b^4))$, if $a$ and $b$ are both odd.\\
The demand that $U_6(P,Q)$ be square is that $P(P^2-Q)(P^2-3Q)=\Box$, which leads to
one of seven cases: $P=a^2$, $P^2-Q=b^2$, with $-2 a^4+3 b^2=\Box$; $P=a^2$, $P^2-Q=-2 b^2$,
with $a^4+3 b^2=\Box$; $P=-a^2$, $P^2-Q=2 b^2$, with $a^4-3 b^2=\Box$; and $P=3 a^2$,
$P^2-Q= \delta b^2$, ($\delta = \pm 1, \pm 2$), with $-\frac{6}{\delta} a^4+b^2=\Box$.
So finitely many parametrizations result (which can easily be obtained if we wish
to do so).\\
The demand that $U_7(P,Q)$ be square is that $P^6-5P^4 Q+6 P^2 Q^2-Q^3 =\Box$,
equivalently, that $1+5 x+6 x^2 +x^3 = y^2$, where $x=-Q/P^2$. This latter elliptic curve
has rank 1, with generator $P_0=(-1,1)$, and trivial torsion. Accordingly,
sequences with $U_7(P,Q)=\Box$ are parametrized by the multiples of $P_0$ on the  above
elliptic curve, corresponding to $(\pm P,Q) = (1,1)$, $(1,5)$, $(2,-1)$, $(5,21)$, $(1,-104)$, 
$(21,545)$, $(52,415)$,...\\
\section{Solution of $U_n(P,Q)=\Box$, $n=8,10,11$} \label{U8 to U11}
Here we shall show the following result:
\begin{u8_result} \label{u8_result}
The only non-degenerate sequences where $\gcd(P,Q)=1$ and $U_n(P,Q)=\Box$, $8 \leq n \leq 11$, 
are given by $U_8(1,-4)=21^2$ and $U_8(4,-17)=620^2$.
\end{u8_result}
\noindent
For each $n \in \{8,10,11\}$ we reduce the solution of our problem to a number of questions, 
all having the following general shape:
\begin{problem} \label{problem}
Let
\begin{equation} \label{gen_xy_elliptic}
\cE\,:\, Y^2+a_1XY+a_3Y=X^3+a_2X^2+a_4X +a_6
\end{equation}
be an elliptic curve defined over $\Q(\alpha)$, where $\alpha$ is a root of a
polynomial $f(X)\in\Z[X]$, irreducible over $\Q$, of degree $d\geq 2$,
and let $\beta,\gamma \in\Q(\alpha)$ be algebraic integers. Find all points
$(X,Y)\in\cE(\Q(\alpha))$ for which $\beta X+\gamma$ is a {\em rational number}.
\end{problem}
\noindent
As mentioned in the Introduction, problems of this type may be attacked with
the Magma routines of Nils Bruin, in particular the routine ``Chabauty" using output from ``PseudoMordellWeilGroup". It is not
{\it a priori} guaranteed that these routines will be successful, particularly over number fields of
high degree, but for our
purposes and the computations encountered in this paper, no
difficulties arose.
\subsection{$n=8$}
The demand that $U_8(P,Q)$ be square is that $P(P^2-2Q)(P^4-4P^2Q+2Q^2)=\Box$. 
In the case $P$ is odd, it follows that $(P,P^2-2Q,P^4-4P^2Q+2Q^2) = (a^2,b^2,c^2)$, 
$(a^2,-b^2,-c^2)$, $(-a^2,b^2,-c^2)$, or $(-a^2,-b^2,c^2)$, where $a$, $b$, $c$ are 
positive integers with $a b$ odd. The latter two possibilities are impossible modulo $4$,
and the first two possibilities lead respectively to:
\begin{eqnarray}
-a^8+2 a^4 b^2 +b^4 & = & 2 c^2 \label{eq1} \\
-a^8-2 a^4 b^2 +b^4 & = & - 2 c^2. \label{eq2}
\end{eqnarray}
We shall see that the only positive solutions to the above equations 
are $(a,b)=(1,1),(1,3)$ and $(1,1)$ respectively, leading
to $(P,Q)=(1,0), (1,-4)$ and $(1,1)$, from which we reject the first one.
The last gives a degenerate sequence.
In the case that $P$ is even, then $Q$ is odd, and $2$ exactly divides both 
$P^4-4P^2Q+2Q^2$ and $P^2-2Q$, forcing $P \equiv 0 \pmod 4$. Put $P=4p$, so that 
$U_8=\Box$ if and only if $p(8p^2-Q)(128p^4-32p^2Q+Q^2)=\Box$, with $\gcd(p,Q)=1$. 
It follows that $(p,8p^2-Q,128p^4-32p^2Q+Q^2)=(a^2,b^2,c^2)$, $(a^2,-b^2,-c^2)$, 
$(-a^2,b^2,-c^2)$, or $(-a^2,-b^2,c^2)$, where $a$, $b$, $c$ are positive integers, 
$(a,b)=1$ and $b c$ is odd. The middle two possibilities are impossible modulo $4$, 
and the remaining two possibilities lead respectively to:
\begin{eqnarray}
-64a^8+16a^4b^2+b^4 & = & c^2 \label{eq3} \\
-64a^8-16a^4b^2+b^4 & = & c^2. \label{eq4}
\end{eqnarray}
We shall see that the only positive solution which leads 
to a desired pair $(P,Q)$ is from the first equation when $(a,b)=(1,5)$, leading to 
$(P,Q)=(4,-17)$.\\ \\
%We reduce the solution of equations (\ref{eq1})-(\ref{eq4}) to the
%solution of a number of problems all of which fit the following general shape: 
%\begin{problem} \label{problem}
%Let 
%\begin{equation} \label{gen_xy_elliptic}
%\cE\,:\, Y^2+a_1XY+a_3Y=X^3+a_2X^2+a_4X +a_6
%\end{equation} 
%be an elliptic curve defined over $\Q(\alpha)$, where $\alpha$ is a root of a 
%polynomial $f(X)\in\Z[X]$, irreducible over $\Q$, of degree $d\geq 2$, 
%and let $\beta,\gamma \in\Q(\alpha)$ be algebraic integers. Find all points 
%$(X,Y)\in\cE(\Q(\alpha))$ for which $\beta X+\gamma$ is a {\em rational number}.
%\end{problem}
%In fact, each equation (\ref{eq1})-(\ref{eq4}) leads to four elliptic curves, 
%and four instances of Problem \ref{problem}.
We shall work in the number field $K=\Q(\phi)$ where $\phi$ is a root of 
$f_\phi(x)=x^4+2 x^2-1$.
The class number of $K$ is 1, the maximal order $\mathcal{O}$ of $K$
is $\Z[\phi]$, and fundamental units of $\mathcal{O}$ are
$\eta_1=\phi$, $\eta_2=2-3\phi+\phi^2-\phi^3$.
The factorization of 2 is $2=\eta_1^{-4} \eta_2^2 (1+\phi)^4$.
\subsubsection{Equation (\ref{eq1})} \label{curve1}
The factorization of (\ref{eq1}) over $K$ is
\[ (b-\phi a^2)(b+\phi a^2)(b^2+(2+\phi^2) a^4) = 2 \Box, \]
and it is easy to see that each of the first two factors is exactly
divisible by $1+\phi$ and the third factor is exactly divisible by $(1+\phi)^2$. Thus the 
gcd of any two (ideal) terms on the left hand side is equal to $(1+\phi)$, and
\begin{equation}
\label{ijeq}
b+\phi a^2 = \pm\eta_1^{i_1}\eta_2^{i_2}(1+\phi)\Box, \qquad b^2+(2+\phi^2)a^4 = 
\pm\eta_1^{j_1}\eta_2^{j_2}(1+\phi)^2\Box,
\end{equation}
where the exponents of the units are 0,1. Specializing $\phi$ at the real root $0.643594...$ 
of $f_\phi(x)$, and using $b>0$, then necessarily the sign on the right hand side must be 
positive.
Taking norms in the first equation gives $2 c^2 = (-1)^{i_1} 2\Box$, so that $i_1=0$.
Applying the automorphism of $K$ defined by $\phi \rightarrow -\phi$, it follows
that $b-\phi a^2 = \eta_2^{-i_2}(1-\phi)\Box = \eta_2^{-i_2+1}(1+\phi)\Box$. Multiplying
this equation by the two displayed equations at (\ref{ijeq}) gives 
$\eta_1^{j_1}\eta_2^{j_2+1}\Box = 2 c^2 = \Box$, so that $j_1=0$, $j_2=1$. We now have
\[ (b+\phi a^2)(b^2+(2+\phi^2)a^4) = \eta_2^j(1+\phi)\Box, \]
with $j=0,1$. Putting $b/a^2 = \delta^{-1} x/(1+\phi)$, where $\delta \in \{1,\eta_2\}$, 
our problem reduces to finding all $K$-points $(x,y)$ on the curves
\[ (x+\phi (1+\phi) \delta)(x^2+(2+\phi^2)(1+\phi)^2 \delta^2) = y^2, \]
subject to $\delta^{-1}x/(1+\phi) \in \Q$. 
For both curves, the Magma routines of Bruin show that solutions to (\ref{eq1})
occur only for $(\pm a, \pm b)=(1,1),(1,3)$.
\subsubsection{Equation (\ref{eq2})} \label{curve2}
As above, (\ref{eq2}) leads to equations
\[ b+\phi^{-1} a^2 = \eta_1^{i_1}\eta_2^{i_2}(1+\phi)\Box, \qquad  b^2+(-2+\phi^{-2}) a^4 = 
\eta_1^{j_1}\eta_2^{j_2}(1+\phi)^2\Box, \]
with exponents $0,1$. Taking norms in the first equation gives $-2 c^2 =(-1)^{i_1} 2 \Box$, so 
that $i_1=1$; and specializing $\phi$ at the real root $-0.643594...$ of $f_\phi(x)$
in the second equation gives $j_1=0$. We thus obtain
\[ (b+\phi^{-1} a^2) (b^2 +(-2+\phi^{-2}) a^4) = \eta_1\eta_2^j(1+\phi)\Box, \]
with $j=0,1$.
Putting $b/a^2 = \delta^{-1} x/(1+\phi)$, where 
$\delta=\eta_1\eta_2^j$, we thus have to find all $K$-points $(x,y)$ on the curves
\[(x+\frac{1+\phi}{\phi} \delta)(x^2+(-2+\frac{1}{\phi^2})(1+\phi)^2 \delta^2) = y^2, \]
such that $\delta^{-1} \frac{x}{1+\phi} \in \Q$, for
$\delta=\eta_1$, $\eta_1 \eta_2$.
The Magma routines show that solutions to (\ref{eq2}) occur only for $(\pm a, \pm b)=(1,1)$.
\subsubsection{Equation (\ref{eq3})} \label{curve3}
As above, (\ref{eq3}) leads to equations
\[ b+2(\phi^3+\phi)a^2 = \eta_1^{i_1}\eta_2^{i_2}\Box, \qquad b^2+8(2+\phi^2)a^4 = \eta_1^{j_1}
\eta_2^{j_2}\Box, \]
with exponents $0,1$. Arguing just as for equation (\ref{eq1}), we deduce $i_1=0$, and $j_1=0$,
$j_2=0$. Thus
\[ (b+2(\phi^3+\phi)a^2) (b^2+8(2+\phi^2)a^4) = \eta_2^j\Box, \]
with $j=0,1$.  For $a \neq 0$, put $b/a^2 = \delta^{-1} x$,
where $\delta=\eta_2^j$, which leads to
seeking all $K$-points $(x,y)$ on the curves
\[ (x+2(\phi^3+\phi)\delta) (x^2 +8(\phi^2+2)\delta^2) = y^2, \]
subject to $\delta^{-1}x \in \Q$, with $\delta = 1$, $\eta_2$.
Magma routines show the only
solutions of (\ref{eq3}) are given by $(\pm a, \pm b) = (0,1),(1,2),(1,5)$,
of which only the third provides a desired pair $(P,Q)$.
\subsubsection{Equation (\ref{eq4})} \label{curve4}
Arguing as in previous cases, we deduce an equation
\[ (b+2(\phi^3+3\phi)a^2) (b^2 +8\phi^2 a^4) = \eta_2^j \Box, \]
with $j = 0,1$.
For $a \neq 0$, put $b/a^2 = \delta^{-1} x$, where $\delta=\eta_2^j$.
This leads to finding all $K$-points $(x,y)$ on the curves 
\[ (x+2(\phi^3+3\phi)\delta) (x^2 +8\phi^2 \delta^2) = \Box, \]
with
$\delta^{-1} x \in \Q$, and $\delta=1$, $\eta_2$. 
Magma routines show the only solution of (\ref{eq4})
occurs for $a=0$, with no solution for $(P,Q)$.
\subsection{$n=10$}
The equation $U_{10}=\Box$ is given by
\begin{equation}
\label{u10PQ}
P (P^4-3 P^2 Q+Q^2) (P^4-5 P^2 Q+5 Q^2) = \Box.
\end{equation}
Our aim is to show that the only integer
solutions are given by $(P,Q)=(1,0)$, $(0,1)$, $(-1,1)$.
The assumption $\gcd(P,Q)=1$ implies that $P$ is coprime
to the second factor on the left at (\ref{u10PQ}), and
may only have the divisor 5 in common with the third factor. Further,
the second and third factors can only have common divisor 2, impossible
for $\gcd(P,Q)=1$.
Putting $(x,y)=(P^2,Q)$, then factorization over $\Z$ implies
\[ x^2-3 x y+y^2 = d_1 z^2, \qquad x^2-5 x y+5 y^2 = d_2 w^2, \]
with $d_1=\pm 1$, $d_2=\pm 1, \pm 5$, giving 8 curves
of genus 1. Most of the possibilities for $(d_1,d_2)$ are
readily eliminated, either by local consideration, or by leading to
rank 0 elliptic curves; the case $(d_1,d_2)=(-1,-5)$ however
resists elementary treatment, apparently leading to a curve of genus 3 
with Jacobian of rank 4. We have found 
it preferable to invoke factorization of (\ref{u10PQ}) over 
$K = \Q(\sqrt{5})$, when (\ref{u10PQ}) becomes
\begin{equation}
P (P^2-\ep^2 Q) (P^2-\ep^{-2} Q) (P^2-\sqrt{5} \ep Q) (P^2 -\sqrt{5} \ep^{-1} Q) = \Box,
\end{equation}
where $\ep=(1+\sqrt{5})/2$ is a fundamental unit of the ring of
integers $\mathcal{O}_K$ of $K$, of class number 1. Denote with a bar 
conjugation under $\sqrt{5} \rightarrow -\sqrt{5}$, so 
that $\bar{\ep} = -\ep^{-1}$.
We have two cases to consider: (1) $(P,5)=1$, and (2) $(P,5)=5$.\\ \\
Case (1): $(P,5)=1$. It follows that $P=\pm p^2$, and
\begin{center}
$\begin{array}{ccccccc}
p^4 - \ep^2 Q & = & \la_1 \al_1^2, & \qquad & p^4 - \ep \sqrt{5} Q & = & \la_2 \al_2^2, \\
p^4 - \ep^{-2} Q & = & \bar{\la_1} \bar{\al_1}^2, & \qquad & p^4 - \ep^{-1} \sqrt{5} Q & = & 
\bar{\la_2} \bar{\al_2}^2 
\end{array}$
\end{center}
where $\la_i$, $\al_i \in \mathcal{O}_K$, with $\la_i$ units.
Equivalently, since $p \neq 0$,
\begin{eqnarray}
1 - \ep^2 q  =  \la_1 \be_1^2, & \qquad & 1 - \ep \sqrt{5} q  =  \la_2 \be_2^2, \nonumber \\
1 - \ep^{-2} q  =  \bar{\la_1} \bar{\be_1}^2, & \qquad & 1 - \ep^{-1} \sqrt{5} q  =  
\bar{\la_2} \bar{\be_2}^2, \label{lambda}
\end{eqnarray}
where $q = Q/p^4 \in {\bf Q}$, $\be_i \in K$, and without loss of generality,
$\la_i \in \{\pm 1, \pm \ep\}$. From the first three of these
equations,
\[ (q-\ep^{-2})(q-\ep^2)(q-\frac{\ep^{-1}}{\sqrt{5}}) = -\la_1 \bar{\la_1} \la_2 
\frac{\ep^{-1}}{\sqrt{5}} \Box = v \frac{\ep^{-1}}{\sqrt{5}} \Box, \]
where $v= -\la_1 \bar{\la_1} \la_2 =\pm \la_2$. Putting
\[ x = \delta q, \quad \delta = v \ep \sqrt{5}, \qquad \mbox{ where} \quad \delta^{-1} x \in 
\Q, \]
then 
\[ (x-v \ep^{-1}\sqrt{5}) (x-v \ep^3\sqrt{5}) (x-v) = \Box. \]
If $v=\pm \la_2 = \pm 1$, then one of the following equations holds:
\begin{eqnarray}
\label{e1}
y^2 & = & (x-(3-\ep)) (x-(3+4\ep)) (x-1) \\
\label{e2}
y^2 & = & (x+(3-\ep)) (x+(3+4\ep)) (x+1)
\end{eqnarray}
under the condition $\frac{\ep^{-1}}{\sqrt{5}} x \in \Q$. Equation (\ref{e1})
defines an elliptic curve of $K$-rank 0, with no corresponding value of $Q$;
equation (\ref{e2}) defines an elliptic curve of positive $K$-rank, and the 
Magma routines show that the only points satisfying the rationality 
condition are $(x,\pm y)=(0,\ep\sqrt{5})$, $(-2-\ep, 1+3\ep)$, corresponding
to $Q=0$ and $Q=1$. \\
If $\la_2=\pm \ep$, then one of the following equations holds:
\begin{eqnarray*}
y^2 & = & (x-\sqrt{5}) (x-\ep^4 \sqrt{5})) (x-\ep) \\
y^2 & = & (x+\sqrt{5}) (x+\ep^4 \sqrt{5})) (x+\ep)
\end{eqnarray*}
under the condition $\frac{\ep^{-2}}{\sqrt{5}} x \in \Q$. Both these curves
have $K$-rank 0, and no solution to our problem arises.\\ \\
Case (2): $(P,5)=5$. We have $P=\pm 5 p^2$, and
\[ (25 p^4-\ep^2 Q)(25 p^4-\ep^{-2} Q)(5 \sqrt{5} p^4 -\ep Q) (5 \sqrt{5} p^4 -\ep^{-1} Q) = 
\pm \Box, \]
so that either $p=0$ (returning the solution $(P,Q)=(0,1)$ to (\ref{u10PQ})), or
\begin{eqnarray}
5 - \ep^2 r = \mu_1 \ga_1^2, & \qquad & \sqrt{5} - \ep r = \mu_2 \ga_2^2, \nonumber \\
5 - \ep^{-2} r = \bar{\mu_1} \bar{\ga_1}^2, & \qquad & \sqrt{5} - \ep^{-1} r = -\bar{\mu_2} 
\bar{\ga_2}^2, \label{mu}
\end{eqnarray}
where $r=Q/(5 p^4) \in {\bf Q}$, $\ga_i \in K$, and $\mu_i$ units of 
$\mathcal{O}_K$, without loss of generality in the set $\{\pm 1, \pm \ep\}$.
From the first three of these equations,
\[ (r-5\ep^{-2}) (r-5\ep^2) (r-\sqrt{5} \ep^{-1}) = -\mu_1 \bar{\mu_1} \mu_2 \ep^{-1} \Box=w 
\ep^{-1} \Box, \]
where $w=-\mu_1 \bar{\mu_1} \mu_2=\pm \mu_2$.  Putting 
\[ x = \eta r, \quad \eta = w \ep, \qquad \mbox{ where} \quad \eta^{-1} x \in \Q, \]
then
\[ (x-5 w \ep^{-1}) (x-5 w \ep^3) (x-w \sqrt{5}) = \Box. \]
If $w=\pm \mu_2 = \pm 1$, then one of the following equations holds:
\begin{eqnarray}
\label{e3}
y^2 & = & (x-5\ep^{-1}) (x-5\ep^3) (x-\sqrt{5}) \\
\label{e4}
y^2 & = & (x+5\ep^{-1}) (x+5\ep^3) (x+\sqrt{5})
\end{eqnarray}
under the condition $\ep^{-1} x \in \Q$.\\
If $w=\pm \ep$, then one of the following equations holds:
\begin{eqnarray}
\label{e5}
y^2 & = & (x-5) (x-5\ep^4) (x-\ep \sqrt{5}) \\
\label{e6}
y^2 & = & (x+5) (x+5\ep^4) (x+\ep \sqrt{5})
\end{eqnarray}
under the condition $\ep^{-2} x \in \Q$. These last four equations
define elliptic curves, the first three of which have positive $K$-rank,
the fourth having $K$-rank 0 (giving no solution to our problem). Magma 
computations show that no solutions arise from the first three curves.
(Computations do disclose the point $(x,y)=(-2\ep, \ep)$ on (\ref{e4})
satisfying the rationality condition, but this point leads to $(P,Q)=(5,10)$,
disallowed for our original problem).
\subsection{$n=11$}
The equation $U_{11}=\Box$ is given by
\begin{equation}
\label{u11PQ}
U_{11}(P,Q)=P^{10}-9 P^8 Q+28 P^6 Q^2-35 P^4 Q^3+15 P^2 Q^4-Q^5 = M^2.
\end{equation}
There is the trivial solution given by $(P,Q)=(0,-1)$, and henceforth we
assume $P \neq 0$.  Our aim is to show that when $\gcd(P,Q)=1$, the only integer
solution of (\ref{u11PQ}) is given by $(P^2,Q)=(1,0)$.
On putting $x=Q/P^2$, $y=M/P^5$, equation (\ref{u11PQ}) becomes
that of a genus 2 curve 
\[ C: y^2 =-x^5+15 x^4-35 x^3+28 x^2 -9 x+1, \]
and Magma computations show that the Jacobian $J$ of $C$ has rank 1,
so that a Chabauty argument can be applied to determine all rational
points of $C$.  A generator is found to be 
$\{ \left( (3+\sqrt{5})/2, (11+5\sqrt{5})/2 \right) - \infty\}$ (where 
$\infty$ is the unique point at infinity on $C$), and Magma tells us that
there is at most one pair of rational points on $C$, which is accordingly
$(0,\pm 1)$ as required. Full details of such an argument may be found
in Flynn, Poonen, Schaefer~\cite{FPS}. For greater transparency, we
outline the method that reduces as in previous cases to finding points
on elliptic curves over a number field under a certain rationality condition.
We shall be working over a Galois field, and this approach requires 
only application of the formal group
of an elliptic curve, conceptually more straightforward than invoking the
formal group of a curve of genus 2.  It is easy to see that
\begin{equation}
\label{Lnorm}
M^2 = u_{11}(P,Q) = \prod_{j=1}^{j=5} ((\zeta^{j/2}+\zeta^{-j/2})^{-2} P^2-Q),
\end{equation}
and thus $U_{11}(P,Q)$ splits completely over the 
real subfield of the cyclotomic field $\Q(\zeta)=\Q(\zeta_{11})$.
Let $\th=\zeta+\zeta^{-1}$, so that
\[ f_\th(\th) = \th^5+\th^4-4\th^3-3\th^2+3\th+1=0, \]
and work in the number field $K=\Q(\th)$, with ring of integers
$\mathcal{O}_K = \Z[\th]$, discriminant $11^4$, and class-number 1.
The unit group is of rank 4, and generators may be taken as:
\[ \ep_1 = -\th, \quad \ep_2 = -\th^2+2, \quad \ep_3 = -\th^4+4\th^2-2, \quad \ep_4 = -\th^3+3
\th, \]
with norms all equal $+1$.  Rewrite (\ref{Lnorm}) in the form
\begin{equation}
\label{MPQeq}
M^2 = \prod_{j=1}^{j=5} L_j(P,Q),
\end{equation}
where $L_j(P,Q)=\th_j P^2-Q$, with $\th_j=(\zeta^{j/2}+\zeta^{-j/2})^{-2}$.
The Galois group of $\Q(\th)/\Q$ is cyclic, generated by the automorphism
$\sigma: \th \rightarrow \th^2-2$, which acts cyclically on the $\th_i$,
and satisfies $\ep_i^\sigma=\ep_{i+1}$ for $i=1,2,3$, $\ep_4^\sigma=\ep_1^{-1}\ep_2^{-1}
\ep_3^{-1}\ep_4^{-1}$.\\
Since $U_{11}(x,1) \equiv 0 \bmod 11^2$ has no solution, 
$M \not \equiv 0 \bmod 11$; so each factor $L_j(P,Q)$ is prime to
11 in $\mathcal{O}_K$. Further, for $j \neq k$, the $K/\Q$ norm of
$\th_j-\th_k$ is $\pm 11$, and it follows that the factors $L_j(P,Q)$
for $j=1,...,5$ are coprime in $\mathcal{O}_K$. Accordingly,
\begin{equation}
\label{L1eq}
L_1(P,Q) = (\th^4-\th^3-2\th^2+\th+1)P^2-Q = (-1)^{i_0} \ep_1^{i_1} \ep_2^{i_2} \ep_3^{i_3} 
\ep_4^{i_4} \Box,
\end{equation}
where $i_0,i_1,...,i_4 \in \{0,1\}$. Since the norm of $L_1(P,Q)$ equals $M^2$,
we must have $i_0=0$.
Now let $^{*}: \Q(\th) \hookrightarrow \R$ be the embedding that sends
$\th$ to the smallest real root $-1.9189859...$ of $f_\th(x)$. It is 
straightforward
to check that $\th_1^* > \th_4^* > \th_5^* > \th_3^* > \th_2^*$, whence
\begin{equation} 
\label{Lineq}
L_1^* > L_4^* > L_5^* > L_3^* > L_2^*.
\end{equation}
Since (\ref{MPQeq}) implies that the exact number of the $L_j^*$ 
that are negative must be even, then (\ref{Lineq}) gives $L_1^* > 0$.
But $\ep_1^* > 0$, $\ep_2^* < 0$, $\ep_3^* < 0$, $\ep_4^* > 0$,
and so necessarily $i_2+i_3$ is even.\\
By applying $\sigma$ repeatedly to (\ref{L1eq}) we obtain:
%\begin{eqnarray*}
%L_1(P,Q) & = & \ep_1^{i_1} \ep_2^{i_2} \ep_3^{i_3} \ep_4^{i_4} \Box \\
%L_2(P,Q) & = & \ep_1^{-i_4} \ep_2^{i_1-i_4} \ep_3^{i_2-i_4} \ep_4^{i_3-i_4} \Box \\
%L_3(P,Q) & = & \ep_1^{-i_3+i_4} \ep_2^{-i_3} \ep_3^{i_1-i_3} \ep_4^{i_2-i_3} \Box \\
%L_4(P,Q) & = & \ep_1^{-i_2+i_3} \ep_2^{-i_2+i_4} \ep_3^{-i_2} \ep_4^{i_1-i_2} \Box \\
%L_5(P,Q) & = & \ep_1^{-i_1+i_2} \ep_2^{-i_1+i_3} \ep_3^{-i_1+i_4} \ep_4^{-i_1} \Box
%\end{eqnarray*}
%so that
\[ \mbox{sgn}(L_2^*)=(-1)^{i_1+i_2}, \; \mbox{sgn}(L_3^*)=(-1)^{i_1}, \; \mbox{sgn}
(L_4^*)=(-1)^{i_4}, \; \mbox{sgn}(L_5^*)=(-1)^{i_3+i_4}. \] 
Listing the 8 possibilities for $(i_1,i_2,i_3,i_4)$, together with the 
corresponding signs of the $L_i$,
%\begin{center}
%$\begin{array}{cccc|ccccc}
%i_1 & i_2 & i_3 & i_4 & \mbox{sgn}L_1^* & \mbox{sgn}L_4^* & \mbox{sgn}L_5^* & \mbox{sgn}L_3^* & 
%\mbox{sgn}L_2^*   \\ \hline
%0 & 0 & 0 & 0 & + & + & + & + & + \\
%0 & 0 & 0 & 1 & + & - & - & + & + \\
%1 & 0 & 0 & 0 & + & + & + & - & - \\
%1 & 0 & 0 & 1 & + & - & - & - & - \\
%0 & 1 & 1 & 0 & + & + & - & + & - \\
%0 & 1 & 1 & 1 & + & - & + & + & - \\
%1 & 1 & 1 & 0 & + & + & - & - & + \\
%1 & 1 & 1 & 1 & + & - & + & - & + 
%\end{array}$
%\end{center}
the only ones that respect the linear ordering (\ref{Lineq}) 
are $(0,0,0,0)$, $(1,0,0,0)$, and $(1,0,0,1)$. The conclusion is that
\[ L_1(P,Q)=\th_1 P^2-Q = \eta \Box, \qquad \eta \in \{1,\;\;\ep_1,\;\;\ep_1\ep_4\}. \]
Certainly
\[ L_1 = \eta \Box, \;\; L_2 = \eta^\sigma \Box, \;\; L_3 = \eta^{\sigma^2} \Box, \;\; L_4 = 
\eta^{\sigma^3} \Box, \;\; L_5 = \eta^{\sigma^4} \Box, \]
and we have various possibilities for producing an elliptic curve cover
of our original equation. 
When $\eta=1$, then
\[ (\th_1 P^2-Q) (\th_2 P^2-Q) (\th_3 P^2-Q) = \Box, \]
so that $x=-Q/P^2$ is the $x$-coordinate of a $K$-rational point
on the elliptic curve
\[ (x+\th_1)(x+\th_2)(x+\th_3) = y^2, \]
with $x \in \Q$. Magma computations show 
that the only such points on this curve are at infinity 
(corresponding to $P=0$), and $(x,y)=(0,\ep_3\ep_4)$,
corresponding to $(P^2,Q)=(1,0)$.\\
When $\eta=\ep_1$ or $\ep_1\ep_4$, then
\[ (\th_1 P^2-Q) (\th_2 P^2-Q) (\th_3 P^2-Q) = \delta \Box, \]
with $\delta=\ep_1\ep_2\ep_3$ or $\ep_1$, respectively; and
$x=-\delta Q/P^2$ is the $x$-coordinate of a $K$-rational point
on the elliptic curve
\[ (x+\delta\th_1)(x+\delta\th_2)(x+\delta\th_3)=y^2, \]
satisfying $\delta^{-1} x \in \Q$. 
These two curves both have positive $K$-rank, and by the Magma routines lead to no non-trivial
solutions for $P$, $Q$. \\ \\
{\bf Acknowledgement}: We thank the anonymous referee for several
useful comments.

\end{document}